\def\C{\mathbb C}

  \def\l{\lambda}

  \def\la{\langle}
  \def\ra{\rangle}

  \def\no{\noindent}
  \def\pf{{\it Proof. }$\;\;$}
  \def\hal{\unskip\nobreak\hfil\penalty50\hskip10pt\hbox{}\nobreak
  \hfill\vrule height 5pt width 6pt depth 1pt\par\vskip 2mm}

    \documentclass[11pt,a4paper]{article}
            \usepackage{amsfonts}

\newtheorem{theorem}{Theorem}
  \newtheorem{thm}{Theorem}[section]
  \newtheorem{propn}[theorem]{Proposition}
  \newtheorem{prop}[thm]{Proposition}
  \newtheorem{lem}[thm]{Lemma}
  
 \newtheorem{coroll}[theorem]{Corollary}

\begin{document}

  \title{Finite subgroups of simple algebraic groups with irreducible centralizers}

  \author{ Martin W.\ Liebeck and Adam R. Thomas}

  \maketitle

\begin{abstract}
We determine all finite subgroups of simple algebraic groups that have irreducible centralizers -- that is, centralizers whose connected component does not lie in a parabolic subgroup.
\end{abstract}

\section{Introduction}

Let $G$ be a simple algebraic group over an algebraically closed field. Following Serre \cite{ser04}, a subgroup of $G$ is said to be $G$-{\it irreducible} (or just {\it irreducible} if the context is clear) if it is not contained in a proper parabolic subgroup of $G$. Such subgroups necessarily have finite centralizer in $G$ (see \cite[2.1]{LT}). In this paper we address the question of which finite subgroups can arise as such a centralizer.
It turns out (see Corollary \ref{centirr} below) that they form a very restricted collection of soluble groups, together with the alternating and symmetric groups $Alt_5$ and $Sym_5$.

The question is rather straightforward in the case where $G$ is a classical group (see Proposition \ref{class} below). Our  main result covers the case where $G$ is of exceptional type.

\begin{theorem}\label{excep}
Let $G$ be a simple adjoint algebraic group of exceptional type in characteristic $p\ge 0$, and suppose $F$ is a finite subgroup of $G$ such that $C_G(F)^0$ is $G$-irreducible. Then $|F|$ is not divisible by $p$, and $F$, $C_G(F)^0$ are as in Tables $\ref{e8}-\ref{g2}$ in Section $\ref{tabs}$ at the end of the paper (one $G$-class of subgroups for each line of the tables).
\end{theorem}

\no {\bf Remarks } (1) The notation for the subgroups $F$ and $C_G(F)^0$ is described at the end of this section; the notation for elements of $F$ is defined in Proposition \ref{cents}.
\begin{itemize}
\item[(2)] The theorem covers adjoint types of simple algebraic groups. For other types, the possible finite subgroups $F$ are just preimages of those in the conclusion.
\item[(3)] We also cover the case where $G = {\rm Aut}\,E_6 = E_6.2$ (see Table \ref{aute6} and Section \ref{e6pf}).
\item[(4)] A complete determination of all $G$-irreducible connected subgroups is carried out in \cite{Tho}.  
\end{itemize}

Every finite subgroup $F$ in Theorem \ref{excep} is contained in a maximal such finite subgroup. The list of maximal finite subgroups with irreducible centralizers is recorded in the next result.

\begin{coroll}\label{maxcor}
Let $G$ be a simple adjoint algebraic group of exceptional type, and suppose $F$ is a finite subgroup of $G$ which is maximal subject to the condition that $C_G(F)^0$ is $G$-irreducible. Then $F$, $C_G(F)^0$ are as in Table $\ref{maxtab}$.
\end{coroll}

\begin{table}
\caption{Maximal finite subgroups $F$ with irreducible centralizer} \label{maxtab}
\[
\begin{array}{|c|c|c|c|}
\hline
G & F & p & C_G(F)^0 \\
\hline
E_8 & 2^4 & p\ne 2 & A_1^8   \\
       & Q_8 &  p\ne 2 &B_2^3  \\
     & 2_-^{1+4} & p\ne 2 &B_1^5  \\
 & Dih_6 & p\ne 2,3 &B_4  \\
 & G_{12} &  p\ne 2,3 & \bar{A}_1A_1A_3  \\
 & Sym_4\times 2 &  p\ne 2,3 & \bar{A}_1A_1A_1 \\
 & SL_2(3) &  p\ne 2,3 & \bar{A}_1A_2 \\
 & 3^2.Dih_8 &  p\ne 2,3 & A_1^2  \\
 & Sym_5 & p\ne 2,3,5 & A_1 \\
 & Q_8 & p=3 & \bar{A}_1D_4 \\
 & Dih_8\times 2 & p=3 & \bar{A}_1^2B_1^2B_2 \\
 & 3 & p=2 & A_8 \\
 & 3^2 & p=2 &  A_2^4  \\
 & 5 & p=2 & A_4^2 \\
 & Frob_{20} & p=3 &  B_2  \\
\hline
\hline
E_7 & 2^3 & p\ne 2 & A_1^7   \\
 & Q_8 & p\ne 2 & \bar{A}_1B_1^4  \\
 & Dih_6 & p\ne 2,3 & A_1A_3 \\
 & Alt_4 & p\ne 2,3 & A_2  \\
 & Sym_4 & p\ne 2,3 & \bar{A}_1A_1 \\
 & 2^2 & p=3 & D_4 \\
 & Dih_8 & p=3 & \bar{A}_1B_1^2B_2 \\
 & 3 & p=2 & A_2A_5 \\
\hline
\hline
E_6 & 2 & p\ne 2 &  A_1A_5 \\
 & Dih_{6} & p\ne 2,3 & A_1A_1  \\
 & 3 & p=2 & A_2^3 \\
\hline
\hline
F_4 & 2^3 & p\ne 2 & A_1^4 \\
     & Q_8 & p\ne 2 & B_1^3 \\
 & Sym_4 & p\ne 2,3 & A_1 \\
 & 3 & p=2 & A_2A_2 \\
 & Dih_8 & p=3 & B_1B_2 \\
\hline
\hline
G_2 & Dih_6 & p\ne 2,3 & A_1  \\
 & 2 & p=3 & A_1A_1  \\
 & 3 & p=2 & A_2  \\
\hline
\end{array}
\]
\end{table}

For the classical groups we prove the following.

\begin{propn}\label{class} Let $G$ be a classical simple algebraic group in characteristic $p \ge 0$ with natural module $V$, and suppose $F$ is a finite subgroup of $G$ such that $C_G(F)^0$ is $G$-irreducible. Then $p \neq 2$ and $F$ is an elementary abelian $2$-group. Moreover, $G \ne SL_n$ and the following hold.  
\begin{itemize}
\item[{\rm (i)}] If $G = Sp_{2n}$, then 
\[
C_{G}(F)^0 = \prod_i Sp_{2n_i} = \prod_i  Sp(W_i),
\]
where $\sum n_i = n$ and $W_i$ are the distinct weight spaces of $F$ on $V$.
\item[{\rm (ii)}] If $G = SO_n$, then   
\[
C_{G}(F)^0 = \prod_i SO_{n_i} = \prod_i SO(W_i),
\]
where $n_i\ge 3$ for all $i$, $\sum n_i = n$ or $n-1$, and $W_i$ are weight spaces of $F$.
\end{itemize}
\end{propn}

\no {\bf Remark } In Section \ref{classproof} we prove a version of this result covering finite subgroups of ${\rm Aut}\,G$ for $G$ of classical type.

\begin{coroll}\label{centirr}
Let $G$ be a simple adjoint algebraic group, and let $D$ be a proper connected $G$-irreducible subgroup. Then the finite group $C_G(D)$ is either elementary abelian or isomorphic to a subgroup of one of the following groups: 
\[2_{-}^{1+4}, \ G_{12}, \ Sym_4 \times 2, \ SL_2(3), \ 3^2.Dih_8, \ Sym_5.\]
\end{coroll}

\no {\bf Notation } Throughout the paper we use the following notation for various finite groups:
\[
\begin{array}{ll}
Z_n, \hbox{ or just }n & \hbox{cyclic group of order } n \\
p^s \ (p \hbox{ prime}) & \hbox{elementary abelian group of order } p^s \\
Alt_n, Sym_n & \hbox{alternating and symmetric groups} \\
Dih_{2n} & \hbox{dihedral group of order } 2n \\
4 \circ Dih_8 & \hbox{order 16 central product with centre } Z_4 \\
2^{1+4}_{-} & \hbox{extra-special group of order 32 of minus type} \\
Frob_{20} & \hbox{Frobenius group of order } 20 \\
G_{12} & \hbox{dicyclic group } \la x,y \ | \ x^6 = 1, x^y = x^{-1}, y^2 = x^3 \ra \hbox{ of order } 12 \\
\end{array}
\]
In the tables in Section \ref{tabs}, and also in the text, we shall sometimes use $\bar A_1$ to denote a subgroup $A_1$ of a simple algebraic group $G$ that is generated by long root subgroups; we use the notation $\bar A_2$ similarly. Also, $B_r$ denotes a natural subgroup of type $SO_{2r+1}$ in a group of type $D_n$. 

We use the following notation when describing modules for a semisimple algebraic group $G$. We let $L(G)$ denote the Lie algebra of $G$. If $\l$ is a dominant weight, then $V_G(\l)$ (or simply $\l$) denotes the rational irreducible $G$-module with high weight $\l$. When $G$ is simple the fundamental dominant weights $\l_i$ are ordered with respect to the labelling of the Dynkin diagrams as in \cite[p. 250]{bourbaki}. If $V_1, \dots, V_k$ are $G$-modules, then $V_1 / \dots / V_k$ denotes a module having the same composition factors as $V_1 + \cdots + V_k$. Finally, when $H$ is a subgroup of $G$ and $V$ is a $G$-module we use $V \downarrow H$ for the restriction of $V$ to $H$.  

\section{Preliminaries}

In this section we collect some preliminary results required in the proof of Theorem~\ref{excep}.

\begin{prop} \label{outeraut}
Let $G$ be a simple adjoint algebraic group in characteristic $p$. Suppose $t \in {\rm Aut}\, G \setminus G$ is such that $t$ has prime order and $C_G(t)^0$ is $G$-irreducible. Then $t$, $C_G(t)^0$ are given in Table $\ref{graphcent}$.  

If $G=D_4$ and $t$ has order $3$ with $C_G(t) = A_2$, there is an involutory graph automorphism of $G$ that inverts $t$ and acts as a graph automorphism on $C_G(t)$.
\end{prop}

\pf The first part follows from \cite[Tables~4.3.3, 4.7.1]{GLS} for $p\ne 2$, from \cite[\S 8]{AS} for $G = D_n$, $p=2$, and from \cite[19.9]{AS} for $G = A_n, E_6$, $p=2$.
The last part follows from \cite[2.3.3]{kleidman}. \hal

\begin{table}\label{graphcent}
\caption{Centralizers of graph automorphisms in simple algebraic groups} 
\[
\begin{array}{|c|c|c|}
\hline
G & \hbox{order of }t & C_{G}(t)^0  \\
\hline
A_{2n} & 2 & B_{n}\,(p\ne 2) \\ 
A_{2n-1} & 2 & C_{n} \\
   && D_{n}\,(p\ne 2) \\
D_n & 2 & B_{n-1} \\
 &&  B_k B_{n-k-1} \,(1\le k\le n-2,\,p\ne 2) \\
D_4 & 3 & G_2 \\
 & &  A_2\,(p\ne 3) \\
E_6 & 2 & F_4 \\
      & & C_4\,(p\ne 2) \\
\hline
\end{array}
\]
\end{table}


\begin{prop}\label{cents}
Let $G$ be a simple adjoint algebraic group of exceptional type in characteristic $p$, and let  $x \in G$ be a non-identity element such that $C_G(x)^0$ is $G$-irreducible. Then $x$ and $C_G(x)$ are as in Table $\ref{irrcents}$; we label $x$ according to its order, which is not divisible by $p$.
\end{prop}

\pf First observe that if $p$ divides the order of $x$ then $C_G(x)^0$ is $G$-reducible by \cite[2.5]{BT}. Hence $x$ is a semisimple element and $C_G(x)^0$ is a semisimple subgroup of maximal rank. It follows from \cite[4.5]{MP} that this implies the order of $x$ is equal to one of the coefficients in the expression for the highest root in the root system of $G$; these are at most 6 for $G=E_8$, and at most 4 for the other types. The classes and centralizers of elements of these orders can be found in \cite[3.1, 4.1]{CG} and \cite[3.1]{CW}. \hal

\begin{table}
\caption{Elements of exceptional groups with irreducible centralizers} \label{irrcents}
\[
\begin{array}{|c|c|c|}
\hline
G & x & C_G(x)  \\
\hline
E_8 & 2A &  A_1E_7 \\
  & 2B & D_8 \\
  & 3A & A_8 \\
  & 3B & A_2E_6 \\
 & 4A & A_1A_7 \\
  & 4B & A_3D_5 \\
  & 5A & A_4^2 \\
  & 6A & A_1A_2A_5 \\
\hline
E_7 & 2A & A_1D_6 \\
       & 2B & A_7.2 \\
       & 3A & A_2A_5 \\
       & 4A & A_1A_3^2.2 \\
\hline
E_6 & 2A & A_1A_5 \\
      & 3A & A_2^3.3 \\
\hline 
F_4 & 2A & B_4 \\
       & 2B & A_1C_3 \\
       & 3A & A_2A_2 \\
       & 4A & A_1A_3 \\
\hline
G_2 & 2A & A_1A_1 \\
       & 3A & A_2 \\
\hline
\end{array}
\]
\end{table}

We shall need a similar result for the group ${\rm Aut}\,E_6 = E_6.2$.

\begin{prop}\label{aute6cent}
Let $G = {\rm Aut}\,E_6 $ and let  $x \in G\setminus G'$ be such that $C_G(x)^0$ is $G'$-irreducible. Then $x$ and $C_{G'}(x)$ are as in Table $\ref{ae6}$.
\end{prop}

\pf If $x$ is an involution then the result follows from Proposition \ref{outeraut}. Suppose now that $1 \ne x^2 \in G'$. Then $x^2$ has order $2$ or $3$ and $C_{G'}(x^2) = A_1 A_5$ or $A_2^3.3$, respectively, by Proposition \ref{cents}. In the former case $x$ acts as a graph automorphism on the $A_5$ factor and $C_{G'}(x) = A_1 C_3$ or $A_1 A_3$ by \cite[Table~4.3.1]{GLS}. Here $A_1 C_3$ is not possible since this lies in a subgroup $F_4$ and hence centralizes an involution in $G \setminus G'$. 

In the case where $C_{G'}(x^2) = A_2^3.3$, the element $x^2$ is of order $3$ in $C_{G'}(x^3)$. By Proposition \ref{cents}, the latter group must be $F_4$ and $C_{F_4}(x^2) = A_2 A_2$.  \hal

\begin{table}
\caption{Elements with irreducible centralizers in $G = {\rm Aut}\,E_6 $} \label{ae6}
\[
\begin{array}{|c|c|}
\hline
x & C_{G'}(x)  \\
\hline
2B & F_4 \\
2C & C_4\,(p\ne 2) \\
4A & A_1A_3\,(p\ne 2) \\
6A & A_2A_2\,(p\ne 3) \\
\hline
\end{array}
\]
\end{table}

We also require information on normalizers of certain maximal rank subgroups. The following proposition can be deduced from \cite[Tables 7--11]{Car} and direct calculation in the Weyl groups of exceptional algebraic groups; many of the results can be found in \cite[Chapter 11]{LSbk}. 

\begin{prop} \label{normalizers}
Let $G$ be a simple algebraic group of exceptional type. Then Table $\ref{norms}$ gives the groups $N_G(M) / M$ (or $N_{{\rm Aut}G}(M) / M$ for $G$ of type $E_6$) for the given maximal rank subgroups $M$ of $G$.
\end{prop}

\begin{table}
\caption{Normalizers of maximal rank subgroups of $G$}\label{norms}
\[
\begin{array}{|c|c|c|}
\hline
G & M & N_G(M) / M  \\
\hline
E_8 & A_8 & 2 \\
& A_2 E_6 & 2 \\
& A_1 A_7 & 2 \\
& A_4^2 & 4 \\
& D_4^2 & Sym_3 \times 2 \\
& A_1^4 D_4 & Sym_4 \\
& A_2^4  & GL_2(3) \\
& A_1^8 & AGL_3(2) \\
\hline
E_7 & A_2A_5 & 2  \\
& A_1^7 & GL_3(2) \\
\hline
E_6 & A_2^3 & Sym_3 \times 2 \\
\hline
F_4 & A_2 A_2 & 2 \\
& D_4 & Sym_3 \\
\hline
G_2 & A_2 & 2 \\
\hline
\end{array}
\]
\end{table}

Next we have a result about the Spin group $Spin_n$ in characteristic $p\ne 2$. Recall that the centre of $Spin_n$ is $2^2$ if $n$ is divisible by 4 and is $Z_2$ if $n$ is odd. In the former case the quotients of $Spin_n$ by the three central subgroups of order two are $SO_n$ and the two half-spin groups $HSpin_n$.

\begin{prop} \label{spincents} Let $G$ be $HSpin_n$ (where $4 | n$) or $Spin_n$ ($n$ odd), in characteristic $p \neq 2$. Let $\la t \ra = Z(G)$ (so that $G / \la t \ra = PSO_n$) and suppose $F$ is a finite $2$-subgroup of $G$ containing $t$ such that $C_G(F)^0$ is $G$-irreducible. Then the preimage of $F /\la t \ra$ in $SO_n$ is elementary abelian. Moreover, an element $e \in F$ has order $2$ if and only if its preimage in $SO_n$ has $-1$-eigenspace of dimension divisible by $4$.
\end{prop}

\pf
If the preimage contains an element $e$ of order greater than 2, then  $C_{SO_n}(e)$ has a nontrivial normal torus, and hence $C_G(F)^0$ cannot be irreducible. The assertion in the last sentence is well known. \hal

In the following statement, by a {\it pure} subgroup of $G$ we mean a subgroup all of whose non-identity elements are $G$-conjugate.

\begin{prop}\label{coh}
Let $G = E_8$ in characteristic $p$.
\begin{itemize}
\item[{\rm (i)}] If $p\ne 2$ then $G$ has two conjugacy classes of subgroups $E \cong 2^2$ such that $C_G(E)^0$ is $G$-irreducible, 
and one class of pure subgroups $E \cong 2^3$; these are as follows:
\[
\begin{array}{lll}
\hline
E & \hbox{elements} & C_G(E)^0 \\
\hline
2^2 & 2B^3 & D_4^2 \\
& 2A^2,2B & A_1^2D_6 \\
2^3 & 2B^7 & A_1^8 \\
\hline
\end{array}
\]
Further, $G$ has no pure subgroup $2^4$.
\item[{\rm (ii)}] If $p\ne 3$ then $G$ has one class of subgroups $E \cong 3^2$ such that $C_G(E)^0$ is $G$-irreducible. For this class, $C_G(E) = A_2^4$.
\end{itemize}
\end{prop}

\pf Part (i) follows from \cite[3.7, 3.8]{CG}. For (ii), let $E = \la x,y\ra < G$ with $E\cong 3^2$ and $C_G(E)^0$ irreducible.
Then $C_G(x) \ne A_8$, so Proposition \ref{cents} implies that $C_G(x) = A_2E_6$, and also that $C_{A_2E_6}(y)^0 = A_2^4$, as required.  \hal

The next result is taken from \cite[Lemma 2.2]{LT}. 

\begin{prop} \label{classirred}
Suppose $G$ is a classical simple algebraic group in characteristic $p \neq 2$, with natural module $V$. Let $X$ be a semisimple connected subgroup of $G$. If $X$ is $G$-irreducible then either
\begin{itemize}
\item[{\rm (i)}] $G = A_n$ and $X$ is irreducible on $V$, or
\item[{\rm (ii)}] $G = B_n, C_n$ or $D_n$ and $V \downarrow X = V_1 \perp \ldots \perp V_k$ with the $V_i$ all non-degenerate, irreducible and inequivalent as $X$-modules.  
\end{itemize}
\end{prop}

\section{Proof of Theorem \ref{excep}}

\subsection{The case $G=E_8$}

We now embark on the proof of Theorem \ref{excep} for the case $G=E_8$. Let $F$ be a finite subgroup of $G$ such that $C_G(F)^0$ is $G$-irreducible. Then $C_G(F)^0$ is semisimple (see \cite[2.1]{LT}). Moreover $C_G(E)^0$ is irreducible for all nontrivial subgroups $E$ of $F$. Also $F$ is a $\{2,3,5\}$-group by Proposition \ref{cents}.

\begin{lem}\label{elab}
If $F$ is an elementary abelian $2$-group, then $F$ is as in Table $\ref{e8}$.
\end{lem}

\pf  Suppose $F \cong 2^k$. If $k\le 2$, or if $k=3$ and $F$ is pure, the conclusion follows from Proposition \ref{coh}(i).

Now assume that $k=3$ and $F$ is not pure. By considering the $2^2$ subgroups of $F$, all of which must be as in (i) of Proposition \ref{coh}, we see that one of these, say $\la e_1,e_2\ra$, is $2B$-pure, so that $C_G(e_1,e_2)^0 = D_4^2$. We have $C_G(e_1)/\la e_1\ra \cong PSO_{16}$, and consider the preimage of $F/\la e_1\ra$ in $SO_{16}$. This preimage is elementary abelian by Proposition \ref{spincents}, so can be diagonalised, and we can take $e_2 = (-1^8,1^8)$. Let $e_3$ be a further element of $F$ that is in class $2A$. Then the $-1$-eigenspace of $e_3$ has dimension 4 or 12, and 
so the fact that $C_G(F)^0$ is $G$-irreducible means that we can take $e_3 = (-1^4,1^4,1^8)$, so that $C_{SO_{16}}(F)^0 = SO_4 SO_4 SO_8$, and so $C_G(E)^0 = A_1^4D_4$ as in Table \ref{e8}.

Next suppose $k\ge 4$. Then $F$ is not pure by Proposition \ref{coh}(i), so $F$ contains a subgroup $\la e_1,e_2,e_3\ra \cong 2^3$ as in the previous paragraph. Arguing as above, we can take a further element $e_4$ of $F$ to be 
$(1^8,-1^4,1^4)$, so that $C_G(e_1,\ldots,e_4)^0 = A_1^8$. There is no possible further diagonal involution in $F$ such that $C(F)^0$ is irreducible, so $k=4$. \hal 

In view of the previous lemma, we assume from now on that $F$ is not an elementary abelian 2-group. Hence if $F$ is a 2-group, it has exponent 4 by Proposition \ref{cents}.

\begin{lem}\label{c42}
Suppose $F$ is a $2$-group, and has no element in the class $4B$. Then one of the following holds:
\begin{itemize}
\item[{\rm (i)}] $F \cong Z_4$, generated by an element in class $4A$, and $C_G(F)^0 = A_1A_7$;
\item[{\rm (ii)}] $F \cong Q_8$ with elements $2A,4A^6$, and $C_G(F)^0 = A_1D_4$.
\end{itemize}
In both cases $F$ is as in Table $\ref{e8}$.
\end{lem}

\pf  Let $e\in F$ have order 4. By hypothesis, $e$ is in class $4A$, so $C_G(e) = A_1A_7$.
There is nothing to prove if $F \cong Z_4$, so assume $|F|>4$ 
and pick $f \in N_F(\la e\ra)\backslash \la e\ra$. If $f \in A_1A_7$ then $C_G(e,f)$ has a normal torus, so $e^f = e^{-1}$ and $f$ induces an involutory graph automorphism on $A_7$ (see Proposition \ref{normalizers}). Hence $C_{A_7}(f)^0 = C_4$ or $D_4$ by Proposition \ref{outeraut}. The subgroup $C_4$ lies in a Levi subgroup $E_6$ of $C_G(A_1) = E_7$ (see the proof of \cite[2.15]{CLSS}), so the irreducibility of $C_G(F)^0$ implies that  $C_{A_7}(f)^0 = D_4$, hence $C_G(e,f)^0 = A_1D_4$. 
Also $\la e,f\ra \cong Q_8$, as shown in the proof of \cite[2.15]{CLSS}. Finally, if $N_F(\la e,f\ra) > \la e,f\ra$ then some element of order 4 in $\la e,f\ra$ has centralizer in $F$ of order greater than 4, which we have seen to be impossible above. Hence $F = \la e,f\ra$. \hal

\begin{lem}\label{4cross2}
Suppose $F$ is a $2$-group and has an element $e$ in the class $4B$. 
If $C_F(e) \ne \la e \ra$, then $F$ is as in Table $\ref{e8}$ (one of the entries $4\times 2$, $Dih_8\times 2$, $4\circ Dih_8$, $Q_8\times 2$, $2_-^{1+4}$).
\end{lem}

\pf Assume that $C_F(e) \ne \la e \ra$, and choose an involution $e_1 \in C_F(e) \backslash \la e\ra$. Diagonalising the preimage of $C_F(e^2)/\la e^2\ra$ in $SO_{16}$ (as in the proof of Proposition \ref{elab}), we can take
\[
e = (-1^6, 1^{10}),\;e_1 = (1^6, -1^4, 1^6), 
\]
so that $\la e,e_1\ra \cong Z_4\times Z_2$ and $C_G(e,e_1)^0 = A_1^2A_3^2$, as in the $4\times 2$ entry in Table \ref{e8}. 

Write $E_0 = \la e,e_1\ra$. If there exists $f \in (F\cap A_1^2A_3^2)\backslash E_0$, then $C_G(e,e_1,f)$ has a nontrivial normal torus, which is a contradiction. Hence $F \cap C_G(E_0)^0 = E_0$. 

Suppose $C_F(e) \ne E_0$. If there is no involution in $C_F(e)\backslash E_0$, then $C_F(e)$ contains a subgroup $Z_4\times Z_4$, which is impossible. So let  $e_2$ be an involution in $C_F(e)\backslash E_0$. As $e_2$ does not centralize $e_1$, we can take
\[
e_2 = (1^6, -1,1^3,-1^3,1^3).
\]
Then $C_G(e,e_1,e_2)^0 = A_3B_1^3$ (where each $B_1$ corresponds to a natural subgroup $SO_3$ in $SO_{16}$),  and $\la e,e_1,e_2\ra = \la e \ra \circ \la e_1,e_2\ra \cong 4 \circ Dih_8$, as in Table \ref{e8}. If $F \ne E_1:=\la e,e_1,e_2\ra$, then there is an element $e_3 \in N_F(E_1)\backslash E_1$, and adjusting by an element of $E_1$ we can take 
\[
e_3 = (-1^3,1^3,-1,1^9).
\]
Then $C_G(E_1,e_3)^0 = B_1^5$ and $\la E_1,e_3\ra \cong 2^{1+4}_-$, as in Table \ref{e8}. Finally, there are no possible further elements of $F$ such that $C_G(F)^0$ is irreducible.

Now suppose $C_F(e) = E_0$ and $F\ne E_0$. Pick $f \in N_F(E_0)\backslash E_0$. Then $f$ centralizes $e^2$, so we can diagonalise in the usual way; adjusting by an element of $E_0$ and using the fact that $C_G(E_0,f)$ has no nontrivial normal torus, we can take $f \in \{e_4,e_5,e_6,e_7\}$, where 
\[
\begin{array}{ll}
e_4=(-1^3,1^3,1^4,-1,1^5), & e_5=(-1^3,1^3,1^4,-1^3,1^3), \\
e_6=(-1^3,1^3,-1,1^3,1^6), & e_7=(-1,1^5,1^4,-1^3,1^3).
\end{array}
\]
Let $E_2 = \la E_0,f\ra$. 

If $f = e_4$ then $C_G(E_2)^0 = B_1^2 \bar{A}_1^2 B_2$ and $E_2 = \la e,e_4\ra \times \la e_1\ra \cong Dih_8\times Z_2$, as in Table \ref{e8}. There are no possible further elements of $F$ in this case.

If $f = e_5$ then $C_G(E_2)^0 = \bar{A}_1^2 B_1^4$ and $E_2 = \la e,e_5\ra \times \la e_1\ra \cong Q_8\times Z_2$, as in Table \ref{e8}. Any further element of $F$ would centralize $e^2$, and hence would violate the fact that $C_F(e) =E_0$. 

Finally, if $f = e_6$ or $e_7$, then $E_2$ is $D_8$-conjugate to $\la e,e_1,e_2\ra$ or
$\la e,e_1,e_4\ra$, cases considered previously. \hal  

\begin{lem}\label{2gp}
If $F$ is a $2$-group, then it is as in Table $\ref{e8}$.
\end{lem}

\pf Suppose $F$ is a 2-group. In view of the previous two lemmas, we can assume that $F$ contains an element $e$ in the class $4B$, and that $C_F(e) =\la e\ra$ --  and indeed that $F$ has no subgroup $Z_4\times Z_2$. 
We can also assume that $F \ne \la e\ra$. Hence there exists $f \in F$ such that $e^f = e^{-1}$. 
As usual we can diagonalise $\la e,f\ra$, and hence take $e = (-1^6,1^{10})$ and $f \in \{f_1,f_2,f_3,f_4,f_5\}$, where 
\[
\begin{array}{ll}
f_1=(-1^3,1^3,-1,1^9), & f_2=(-1^3,1^3,-1^3,1^7), \\
f_3=(-1^3,1^3,-1^5,1^5), & f_4=(-1,1^5,-1^3,1^7), \\
f_5 = (-1,1^5,-1^5,1^5). &
\end{array}
\]
Moreover, the fact that $F$ has no subgroup $Z_4\times Z_2$ implies that $F = \la e,f\ra$.

It is easily seen that the possibilities for $F$ and $C_G(F)^0$ are as follows:
\[
\begin{array}{cc}
\hline
F & C_G(F)^0 \\
\hline
\la e,f_1\ra \cong Dih_8 & B_1^2B_4 \\
\la e,f_2\ra \cong Q_8 & B_1^3B_3 \\
\la e,f_3\ra \cong Dih_8 & B_1^2B_2^2 \\
\la e,f_4\ra \cong Dih_8 & B_1B_2B_3 \\
\la e,f_5\ra \cong Q_8 & B_2^3 \\
\hline
\end{array}
\]
All these possibilities are in Table \ref{e8}. \hal 


\begin{lem}\label{3gp}
If $F$ is a $3$-group, then $F = 3$ or $3^2$ is as in Table $\ref{e8}$.
\end{lem}

\pf Suppose $F$ is a 3-group. It has exponent 3, by Proposition \ref{cents}. If $F=3$ or $3^2$, it is as in Table \ref{e8} by Propositions \ref{cents} and \ref{coh}(ii), so assume $|F|>9$. 

If $F$ has an element $e$ with $C_G(e) = A_8$, then there is an element $f \in C_F(e)\setminus\la e\ra$, and $C_{A_8}(f)$ is reducible in $A_8$, a contradiction. Hence all non-identity elements of $F$ have centralizer $A_2E_6$ (see Proposition \ref{cents}). In particular, they have trace 5 on the adjoint module $L(G)$ (see \cite[3.1]{CG}).

Let $V$ be a normal subgroup of $F$ with $V\cong 3^2$. Then $C_G(V) = A_2^4$ by Proposition \ref{coh}(ii). If $f \in F\setminus V$, then $f$ acts as a 3-cycle on the four $A_2$ factors, so $C_G(V,f)^0 = \bar A_2A_2$. On the other hand, 
since every non-identity element has trace 5 on $L(G)$, we have
\[
\dim C_{L(G)}(V,f) = \frac{1}{27}(248+26\cdot 5) = 14.
\]
This is a contradiction, showing that $|F|>9$ is impossible. \hal

\vspace{2mm}
From now on, we assume that $F$ is not a 2-group or a 3-group. Let $J = {\rm Fit}(F)$, the Fitting subgroup of $F$.

\begin{lem}\label{fit2}
Suppose $J$ is a nontrivial $2$-group. Then $F$ is as in Table $\ref{e8}$.
\end{lem}

\pf By Lemma \ref{2gp}, $J$ and $C_G(J)^0$ are as in Table \ref{e8}. Also $C_F(J) \le J$, so $F$ contains an element $x$ of order $r=3$ or 5 acting nontrivially on $J$ and as a graph automorphism of $C_G(J)^0$. By inspection of Table \ref{e8}, the possibilities for $J$ with these properties are as follows: 
\[
\begin{array}{c c c}
\hline
J & C_G(J)^0 & r \\
\hline
2^2 & D_4^2 & 3 \\
2^3 & A_1^4D_4 & 3 \\
    & A_1^8 & 3,5 \\
2^4 & A_1^8 & 3,5 \\
Q_8 & A_1D_4 & 3 \\
       & B_2^3 & 3 \\
       & B_1^3B_3 & 3 \\
4\circ Dih_8 & A_3B_1^3 & 3 \\
Q_8\times 2 & A_1^2B_1^4 & 3 \\
2_-^{1+4} & B_1^5 & 3,5 \\
\hline
\end{array}
\]
For the last five cases, where $C_G(J)^0 \triangleright B_1^r$ (or $B_2^r$), $C_G(J,x)^0$ has a factor $B_1$ (or $B_2$) which is a diagonal subgroup of this, and so $C_G(J,x)^0$ is a reducible subgroup of $D_8$ in these cases, by Proposition \ref{classirred}. Also if $C_G(J)^0 = A_1^8$, then $N_G(A_1^8)/A_1^8 \cong AGL_3(2)$ by Proposition \ref{normalizers}, so $x$ acts as either a product of two 3-cycles or a 5-cycle on the eight $A_1$ factors. We claim that again $C_G(J,x)^0$ is reducible. To see this, regard $A_1^8$ as a subgroup of $D_8$ corresponding to $SO_4^4$ in $SO_{16}$. When $x$ is a 5-cycle, $C_J(x)^0 = \bar A_1^3A_1$, where the last factor is diagonal in $\bar A_1^5$. Since $AGL_3(2)$ is 3-transitive, we can choose these 5 factors $\bar A_1$ so that the diagonal subgroup $A_1$ is contained in $\bar A_1B_1B_1$, hence $C_J(x)^0$ is reducible in $D_8$. And when $x$ has order 3, $C_J(x)^0 = \bar A_1^2A_1^2$ where each of the last two $A_1$ factors is diagonal in $\bar A_1^3$. There are two possible actions of the subgroup $A_1^2 < \bar{A}_1^6 < D_6$ on the $12$-dimensional natural module, namely $(1,1)^3$ or $(2,0) + (1,1) + (0,2) + (0,0)^2$. In both cases the subgroup $A_1^2$ is $D_6$-reducible by Proposition \ref{classirred} and hence $C_J(x)^0$ is $D_8$-reducible. 

This leaves the following possibilities remaining, all with $r=3$:
\[
\begin{array}{c c}
\hline
J & C_G(J)^0 \\
\hline
2^2 & D_4^2  \\
2^3 & A_1^4D_4  \\
Q_8 & A_1D_4  \\
\hline
\end{array}
\]

Suppose $J = 2^2$, $C_G(J)^0 = D_4^2$. Then $x$ induces a triality automorphism on both $D_4$ factors (see Proposition \ref{normalizers}). So by Proposition \ref{outeraut}, $C_G(J,x)^0 = G_2G_2$, $A_2A_2$ or $A_2G_2$. In the first case $G_2G_2 < D_7 < D_8$, so is reducible. Hence if $F = \la J,x\ra$ we have the possibilities
\[
F = \la J,x\ra \cong Alt_4, \; C_G(F)^0 = A_2A_2 \hbox{ or } A_2G_2,
\]
both in Table \ref{e8}. Now assume $F \ne \la J,x\ra$. Then $F = \la J,x,t\ra \cong Sym_4$, where $t$ is an involution inverting $x$ (since ${\rm Fit}(F) \cong 2^2$). If $C_G(J,x)^0 = A_2A_2$ then $t$ acts as a graph automorphism on each $A_2$ factor (see Proposition \ref{outeraut}), and so $C_G(F)^0 = A_1A_1$; and if $C_G(J,x)^0 = A_2G_2$ then $t$ acts as a graph autormorphism on the $A_2$ factor and centralizes the $G_2$ factor, so $C_G(F)^0 = A_1G_2$. Hence we have the possibilities
\[
F = \la J,x,t\ra \cong Sym_4, \; C_G(F)^0 = A_1A_1 \hbox{ or } A_1G_2,
\]
both in Table \ref{e8}.

Next suppose  $J = 2^3$, $C_G(J)^0 = A_1^4D_4$. Then $x$ acts as a 3-cycle on the $A_1$ factors and as a triality on $D_4$, so $C_G(J,x)^0 = \bar A_1A_1G_2$ or $\bar A_1A_1A_2$ (where $\bar A_1$ denotes a fundamental $SL_2$ generated by a root group and its opposite). The first subgroup is reducible as it is contained in a subgroup $D_7$ of $D_8$. So if $F = \la J,x\ra$, we have 
\[
F = \la J,x\ra \cong 2\times Alt_4, \; C_G(F)^0 = \bar A_1A_1A_2,
\]
as in Table \ref{e8}.  Now assume $F \ne \la J,x\ra$, and let $\la v\ra = Z(\la J,x\ra)$. As $F/J$ is isomorphic to a subgroupof  $GL_3(2)$ with no nontrivial normal 2-subgroup, we have $F/J \cong Dih_6$ and $F = \la J,x,t\ra$ where $x^t=x^{-1}$ and $t^2=1$ or $v$. Such an element $t$ centralizes both $A_1$ factors of $C_G(J,x)^0$, and induces a graph automorphism on the $A_2$ factor, so $C_G(F)^0 = \bar A_1A_1A_1$. This is contained in the above centralizer $G_2A_2$ of an $Alt_4$ subgroup, and $\bar A_1A_1$ centralizes an involution in $G_2$. Hence in fact $t^2=1$ and we have 
\[
F = \la J,x,t \ra \cong 2\times Sym_4, \; C_G(F)^0 = \bar A_1A_1A_1,
\]
as in Table \ref{e8}.

Finally, suppose $J = Q_8$, $C_G(J)^0 = A_1D_4 < A_1A_7$. Then $x$ induces triality on the $D_4$ factor (see \cite[2.15]{CLSS}), so $C_G(J,x)^0 = A_1G_2$ or $A_1A_2$. The first subgroup is reducible in $A_1A_7$, so if 
$F = \la J,x\ra$, we have 
\[
F = \la J,x\ra =Q_8.3 \cong SL_2(3), \; C_G(F)^0 = \bar A_1A_2,
\]
as in Table \ref{e8}. If $F \ne \la J,x\ra$ then $F$ has an element $t$ inducing a graph automorphism on the $A_2$ factor, so $C_G(J,x,t)^0 = A_1A_1$, which is reducible in $A_1A_7$. This completes the proof. \hal 

\begin{lem}\label{fit3}
If $|J|_3=3$, then $F$ is as in Table $\ref{e8}$.
\end{lem}

\pf Assume $|J|_3=3$, and let $x \in J$ be of order 3. 

Suppose first that $|F|_3=3$ also. As $F$ has no element of order 15 we have $|F|_5=1$, so $F/\la x\ra$ is a 2-group. The case where $|F| = 3$ is in Table \ref{e8}, so assume $|F|>3$.

Suppose $C_G(x) = A_2E_6$. If $t$ is an involution in $C_F(x)$, then $C_G(x,t) = A_2A_1A_5$; moreover $C_F(x)$ has no element of order 4 (as $F$ has no element of order 12), and no subgroup $V \cong 2^2$ (as $C_G(x,V)$ would have a normal torus). Hence $C_F(x)$ has order 3 or 6, and so $|F|$ is 6
or 12. 

If $|F| = 6$, then either $F\cong Z_6$, $C_G(F) = A_1A_2A_5$, or $F = \la x,t\ra \cong Dih_6$ with $t$ inducing graph automorphisms on both factors of $C_G(x) = A_2E_6$ (see Proposition \ref{normalizers}), in which case $C_G(F)^0 = A_1F_4$ or $A_1C_4$. All these possibilities are in Table \ref{e8}.

If $|F|=12$, then $F = \la y,u\ra$ where $y$ has order 6, $y^u = y^{-1}$ and $u^2=1$ or $y^3$. Then $u$ induces a graph automorphism on the $A_2,A_5$ factors of $C_G(y) = A_1A_2A_5$, so $C_G(F)^0 = \bar A_1A_1C_3$ or 
$\bar A_1A_1A_3$. The subgroup $\bar A_1A_1C_3$ is contained in $A_1F_4$ and in $A_1C_4$, while the subgroup 
$\bar A_1A_1A_3$ is contained in neither. Hence $u$ is an involution in the first case, and has order 4 in the second. This gives the possibilities
\[
F = \la y,u\ra \cong Dih_{12} \hbox{ or }G_{12},\; C_G(F)^0 = \bar A_1A_1C_3 \hbox{ or } \bar A_1A_1A_3 \hbox{ (resp.)},
\]
both in Table \ref{e8}.

Next suppose that $C_G(x) = A_8$. Then $C_F(x) = \la x\ra$, so $F= \la x,t\ra \cong Dih_6$ where $t$ induces a graph automorphism on $A_8$, giving $C_G(F)^0 = B_4$. This completes the case where $|F|_3=3$.

Finally, suppose $|F|_3 = 3^2$, and let $\la x,y\ra$ be a Sylow 3-subgroup of $F$. 
Again, $|F|_5=1$. If $F = \la x,y \ra$ then $C_G(F)^0 = A_2^4$ by Proposition \ref{coh}(ii), as in Table \ref{e8}. Otherwise, as $C_F(J)\le J$ there must be a subgroup $V\cong 2^2$ of $J$ such that $y$ acts nontrivially on $V$. But then $C_G(x,V)^0 = C_{A_2E_6}(V)^0$ is reducible, a contradiction. \hal 

\begin{lem}\label{fit32}
If $|J|_3=3^2$, then $F$ is as in Table $\ref{e8}$.
\end{lem}

\pf Let $V \cong 3^2$ be a Sylow 3-subgroup of $J$ (also of $F$, by Lemma \ref{3gp}). By Propositions \ref{coh}(ii) and \ref{normalizers}, $C_G(V) = A_2^4$ and $N_G(A_2^4)/A_2^4 \cong GL_2(3)$. There is no involution in $C_F(V)$, so $J = {\rm Fit}(F) = V$ and $F/J$ is a nontrivial 2-subgroup of $GL_2(3)$.

Suppose first that $|F/J|=2$. There are two classes of involutions in $GL_2(3)$, with representatives $i=-I$ and $t = \pmatrix{-1&0\cr 0&1}$. Then $i$ induces a graph automorphism on each $A_2$ factor of $C_G(V)$, so $C_G(V,i)^0 = A_1^4$; and $t$ fixes two $A_2$ factors, inducing a graph automorphism on one of them, so $C_G(V,t)^0 = A_1\bar A_2A_2$. Both these groups $F = 3^2.2$ are in Table \ref{e8}.

If $F/J\cong 2^2$, we can take $F = \la V,i,t\ra$ and so $C_G(F)^0 = A_1^2A_1$, as in Table \ref{e8}.

Next suppose $F/J\cong Z_4$. There is one class of elements of order 4 in $GL_2(3)$, with representative $u = \pmatrix{0&1\cr -1&0}$; this swaps two pairs of $A_2$ factors, and squares to $i$. Hence $C_G(F)^0 = C_G(V,u)^0 = A_1^2$, as in Table \ref{e8}. 

If $F/J\cong Dih_8$, we can take $F = \la V,u,t\ra$, and again $C_G(F)^0 = A_1^2$.

Now suppose $F/J\cong Q_8$. Then $F/J$ acts transitively on the four $A_2$ factors, and contains $i$, so $C_G(F)^0 = A_1$, a diagonal subgroup of $A_1^4 < A_2^4$. We claim that $C_G(F)^0$ is reducible. To see this, observe that $A_1^4$ centralizes the involution $i$; this involution corresponds to $w_0$, the longest element of the Weyl group of $G$, and so $C_G(i) = D_8$. Now it is easy to check that $C_G(F)^0 = A_1$ is reducible in this $D_8$. Indeed, $A_1^4$ acts as $(1,1,1,1)$ on the natural module for $D_8$ and hence a diagonal subgroup $A_1$ acts as $4 + 2^3 + 0^2$. Thus $F/J \cong Q_8$ is impossible, and we have now covered all possibilities for $F/J$. \hal

\begin{lem}\label{fit5}
If $|J|_5\ge 5$, then $F$ is as in Table $\ref{e8}$.
\end{lem}

\pf Suppose $|J|_5\ge 5$, and let $x \in J$ have order 5. As $C_G(x) = A_4^2$  by Proposition \ref{cents}, there is no element of order 5 in $C_F(x)\setminus \la x\ra$, and so $\la x \ra$ is a Sylow 5-subgroup of $F$.

As $F$ has no element of order 10 or 15, we have $C_F(x) = \la x\ra$, and $|F| = 10$ or 20. By Proposition  \ref{normalizers}, $N_G(A_4^2)/A_4^2 = \la t \ra\cong Z_4$, where $t$ interchanges the two $A_4$ factors and $t^2$ induces a graph automorphism on both. Hence $F$ is either $Dih_{10}$ or $Frob_{20}$, and $C_G(F)^0 = B_2B_2$ or $B_2$, respectively, as in Table \ref{e8}. \hal

\vspace{2mm}
Lemmas \ref{fit2} -- \ref{fit5} cover all cases where the Fitting subgroup $J$ is nontrivial.

\begin{lem}\label{j1}
Suppose $J = {\rm Fit}(F) = 1$. Then $F = Alt_5$ or $Sym_5$ is as in Table $\ref{e8}$.
\end{lem}

\pf In this case $S:={\rm soc}(F)$ is a direct product of non-abelian simple groups. As $5^2$ does not divide $|F|$, in fact $S$ is simple. Proposition 1.2 of \cite{geomded} shows that $S \cong Alt_5$ or $Alt_6$. 

Suupose $S \cong Alt_5$. Then $S$ has subgroups $D\cong Dih_{10}$ and $A \cong Alt_4$, and by what we have already proved, these subgroups are in Table \ref{e8}. Hence the involutions in $S$ are in the class $2B$ (since those in $A$ are in this class). If the elements of order 3 in $S$ are in class $3A$ (with centralizer $A_8$), then from \cite[3.1]{CG} we see that the traces of the elements in $S$ of orders $2, 3, 5$  on $L(G)$ are $-8,-4,-2$ respectively, and hence
\[
\dim C_{L(G)}(S) = \frac{1}{60}(248-8\cdot 15-4\cdot 20-2\cdot 24) = 0,
\]
which is a contradiction. It follows that the elements of order 3 in $S$ are in the class $3B$, with centralizer $A_2E_6$ and trace 5, so that 
\[
\dim C_{L(G)}(S) = \frac{1}{60}(248-8\cdot 15+5\cdot 20-2\cdot 24) = 3.
\]
Since $C_G(D)^0 = B_2B_2<A_4A_4$, it follows that $C_G(S)^0 = A_1$, embedded diagonally and irreducibly in $A_4A_4$. Also $C_G(A_1) = Sym_5$ by \cite[1.5]{geomded}. Hence $F = Alt_5$ or $Sym_5$ and $C_G(F)^0 = A_1$, as in Table \ref{e8}.

Finally, suppose $S \cong Alt_6$ and choose a subgroup $T<S$ with $T \cong Alt_5$. By the above, $C_G(T)^0 = A_1$ and so $C_G(S)^0$ must also be $A_1$. But as observed before, $C_G(A_1) = Sym_5$, a contradiction. \hal

\vspace{4mm}
We have now established that $F$ and $C_G(F)^0$ must be as in Table \ref{e8}.  To complete the proof of Theorem \ref{excep}, we need to establish that all these examples exist. This is proved in the following lemma.

\begin{lem}\label{exist}
Let $F$ and $C_G(F)^0$ be as in Table $\ref{e8}$. Then $C_G(F)^0$ is $G$-irreducible.
\end{lem}

\pf Any subgroup containing a $G$-irreducible subgroup is itself $G$-irreducible. Thus we need only consider the subgroups $C_G(F)^0$ for which $F$ is maximal. These subgroups are given in Table \ref{maxtab}; let $X$ be such a subgroup $C_G(F)^0$. 

Firstly, if $X$ has maximal rank then $X$ is clearly $G$-irreducible. For the subgroups not of maximal rank we use the fact that a subgroup with no trivial composition factors on $L(G)$ is necessarily $G$-irreducible (since the Lie algebra of the centre of a Levi subgroup gives a trivial composition factor). It thus remains to show $X$ has no trivial composition factors on $L(E_8)$. We find the composition factors of $X$ on $L(G)$ by restriction from a maximal rank overgroup $Y$, as given in the last column of Table \ref{e8}. The restrictions $L(G) \downarrow Y$ are given in \cite[Lemma 11.2, 11.3]{LSbk} for all of the maximal rank overgroups $Y$ except for $A_1 A_7$, $A_1^4 D_4$ and $A_2^4$. The latter subgroups are contained in $A_1 E_7$, $D_4^2$ and $A_2 E_6$, respectively, and it is straightforward to compute their composition factors on $L(G)$. 

We finish the proof with two examples of how to calculate the composition factors of $L(G) \downarrow X$ from those of a maximal rank overgroup $Y$. The others all follow similarly and in each case there are no trivial composition factors. 

For the first example, let $X = B_2^3$ so $p \neq 2$ and $X$ is contained in the maximal rank overgroup $D_8$. From \cite[Lemma 11.2]{LSbk}, 
\[ L(G) \downarrow D_8 = V(\lambda_2) + V(\lambda_7),
\] 
the sum of the exterior square of the natural module for $D_8$ and a spin module. To find the restriction of the spin module $V_{D_8}(\lambda_7)$ to $X$ we consider the chain of subgroups $X < B_2 D_5 < B_2 B_5 < D_8$. By \cite[Lemma 11.15(ii)]{LSbk}, $V_{D_8}(\lambda_7) \downarrow B_2 B_5 = 01 \otimes \lambda_5$. Also, $V_{B_5}(\lambda_5) \downarrow D_5 = \lambda_4 + \lambda_5$ and $V_{D_5}(\lambda_i) \downarrow B_2^2 = 01 \otimes 01$ for $i = 4, 5$. 
Therefore, 
\[L(G) \downarrow X = \bigwedge\nolimits^{\!2} (10 \otimes 00 \otimes 00 + 00 \otimes 10 \otimes 00 + 00 \otimes 00 \otimes 10 + 0) + (01 \otimes 01 \otimes 01)^2\]
and this has no trivial composition factors.

For the second example, let $X = A_1 D_4$. Here $p = 3$ and $X$ is contained in a maximal rank subgroup $A_1 A_7$. Then using the restriction $L(G) \downarrow A_1 E_7$ given in \cite[Lemma 11.2]{LSbk} we find 
\[L(G) \downarrow A_1A_7 = 2\otimes0 + 1\otimes \l_2 + 1\otimes \l_6 + 0\otimes (\lambda_1 + \lambda_7) + 0\otimes \l_4.
\]
It is sufficient to show there are no trivial composition factors for $D_4$ acting on $V_{A_7}(\l)$ for $\l = \l_1 + \l_7$ and $\l_4$. By weight considerations, the first module restricts to $D_4$ as $V(2\l_1)+V(\l_2)$ and the second as $V(2\l_3) + V(2\l_4)$. Hence $L(G) \downarrow X$ has no trivial composition factors. \hal

This completes the proof of Theorem \ref{excep} for $G=E_8$. 

\subsection{The case $G=E_7$}

In this section we prove Theorem \ref{excep} for $G = E_7$, of adjoint type. Let $F$ be a finite subgroup of $G$ such that $C_G(F)^0$ is $G$-irreducible. As before, $C_G(F)^0$ is semisimple and $C_G(E)^0$ is $G$-irreducible for all nontrivial subgroups $E$ of $F$. Also $F$ is a $\{2,3\}$-group by Proposition \ref{cents}. 

\begin{lem} \label{e7lem1}
If $F$ is an elementary abelian $2$-group , then $F$ is as in Table $\ref{e7}$.
\end{lem}

\pf We may suppose that $|F|>2$. 
If $F$ has an element $e$ in the class $2B$, then any further element $f \in F \backslash \la e \ra$ must lie in $C_G(e)\backslash C_G(e)^0 = A_7.2\backslash A_7$, and hence $F = \la e,f \ra \cong 2^2$; moreover $C_G(F)^0 = D_4$, as in the proof of Lemma \ref{c42}. Hence $F$ is as in Table \ref{e7}.

So now suppose that $F$ is $2A$-pure. Let $1\ne e\in F$ and $e_1 \in F \backslash \la e \ra$. Then $C_G(e) = A_1D_6$, and diagonalising in $SO_{12}$ as in Lemma \ref{elab}, we can take $e_1 = (-1^4,1^8)$. Hence $C_G(e,e_1)^0 = A_1^3D_4$. If there is an element $e_2 \in E\backslash \la e,e_1\ra$, then we can take $e_2 = (1^4,-1^4,1^4)$, and so $C_G(e,e_1,e_2)^0 = A_1^7$. Both these possibilities are in Table \ref{e8}, and there are no further possible elements in $F$. \hal

\begin{lem} \label{e7lem2}
If $F$ is a $2$-group containing an element of order $4$, then $F$ is as in Table $\ref{e7}$.
\end{lem}

\pf
Let $e \in F$ of order 4. By Proposition \ref{cents} we have $C_G(e)^0 = A_1A_3^2$. Suppose $F \ne \la e\ra$, so there exists $f \in F$ such that $e^f = e^{-1}$. Now $C_G(e^2) = A_1D_6$, and diagonalising in $SO_{12}$ as in Lemma \ref{elab}, we may take $e = (-1^6,1^6)$ and $f  \in \{f_1,f_2\}$, where
\[
f_1 = (-1,1^5,-1^3,1^3),\;\;f_2 = (-1^3,1^3,-1^3,1^3).
\]
If $f=f_1$ then $C_G(e,f)^0 = \bar{A}_1B_1^2B_2$ and $\la e,f\ra \cong Dih_8$; and if $f=f_2$ then $C_G(e,f)^0 = \bar{A}_1 B_1^4$ and $\la e,f\ra \cong Q_8$. Both possibilities are in Table \ref{e7}. Finally, there are no possible further elements of $F$, as can be seen by diagonalising in the usual way. \hal

In view of the previous two lemmas we assume from this point that $F$ contains an element $x$ of order $3$. Let $J$ be the Fitting subgroup of $F$. Note that $F$ does not contain an element of order $6$ by Proposition \ref{cents}. Therefore $J$ is a $2$-group or a $3$-group. 

\begin{lem} 
If $J$ is a $3$-group then $F$ and $C_G(F)^0$ are as given in Table $\ref{e7}$. 
\end{lem}

\pf
Suppose $|J|=3$. If $|F| = 3$ then by Proposition \ref{cents} we have $C_G(F) = A_2 A_5$. Otherwise $F \cong Dih_6$ and $C_G(F)^0 = A_1 C_3$ or $A_1 A_3$. 

Finally, $|J| > 3$ is impossible because the centralizer of an element of order $3$ in $A_2 A_5$ is not $A_2 A_5$-irreducible. \hal

We may now assume that $J$ is a $2$-group. By Lemmas \ref{e7lem1} and \ref{e7lem2}, $J$ is as in Table \ref{e7} and the action of $x$ shows that the only possibilities are $J \cong 2^2$, $2^3$ or $Q_8$.  

\begin{lem} 
If $J \cong 2^2$ then $F$ and $C_G(F)^0$ are as given in Table $\ref{e7}$. 
\end{lem}

\pf
Suppose $C_G(J)^0 = A_1^3 D_4$. By Proposition \ref{normalizers}, $N_G(A_1^3 D_4) / A_1^3 D_4 \cong Sym_3$ acting simultaneously on both the $A_1^3$ and the $D_4$ factors. Therefore $C_G(J,x)^0 = A_1 A_2$ or $A_1 G_2$ with $\la J, x \ra \cong Alt_4$. The subgroup $A_1 G_2$ is $A_1 D_6$-reducible by Proposition \ref{classirred}, and therefore does not appear in Table \ref{e7}. If $F \ne \la J, x \ra$ then we must have $F \cong Sym_4$ with $C_G(F)^0 = A_1 A_1$. 

Now suppose $C_G(J)^0 = D_4 < A_7$. By \cite[Lemma~2.15]{CLSS}, we have $N_G(D_4) / (D_4 \times C_G(D_4)) \cong Sym_3$. Therefore $C_G(J,x)^0 = A_2$ or $G_2$. The subgroup $G_2$ is $A_7$-reducible and therefore does not appear in Table \ref{e7}. If $F \ne \la J,x \ra$ then $F \cong Sym_4$ with $C_G(F)^0 = A_1$. \hal

\begin{lem}
There are no possible subgroups $F$ with $J \cong 2^3$ or $Q_8$.  
\end{lem}

\pf
Suppose $J \cong 2^3$ so $C_G(J)^0 = A_1^7$. By Proposition \ref{normalizers} we have $N_G(A_1^7) / A_1^7 \cong GL_3(2)$. The element $x \in F$ therefore acts as a product of two disjoint $3$-cycles on the seven $A_1$ factors. But the centralizer $C_G(J,x)^0$ is then $A_1 D_6$-reducible by an argument in the first paragraph of the proof of Lemma \ref{fit2}.

Finally, if $J \cong Q_8$ and $C_G(J)^0 = A_1 B_1^4$, then $C_G(J,x)^0 = A_1 B_1 B_1$ which is clearly $A_1 D_6$-reducible. \hal

\vspace{4mm}
The proof of Theorem \ref{excep} for $G=E_7$ is now complete, apart from showing that all the subgroups $C_G(F)^0$ in Table \ref{e7} are $G$-irreducible. This is proved in similar fashion to Lemma \ref{exist}.

\subsection{The case $G={\rm Aut}\,E_6$} \label{e6pf}

Let $G={\rm Aut}\,E_6 = E_6.2$, and let $F$ be a finite subgroup of $G$ such that $C_{G'}(F)^0$ is $G'$-irreducible. 

\begin{lem} \label{e6order4}
If $F$ has an element $x$ of order $4$, then $F = \langle x \rangle$ and $C_{G'}(F)^0 = A_1 A_3$. 
\end{lem}

\pf
By Proposition \ref{aute6cent}, $C_G(x)^0 = A_1 A_3$. By Proposition  \ref{classirred}, $A_1 A_3$ contains no proper $A_1 A_5$-irreducible connected subgroups and therefore $F = \langle x \rangle$ as claimed. \hal

We now assume that $F$ has no element of order 4. 

\begin{lem}
If $F$ is an elementary abelian $2$-group then it appears in Table $\ref{aute6}$. 
\end{lem}

\pf
If $|F| = 2$ then $F$ and $C_G(F)^0$ are as in Table \ref{aute6} by Proposition \ref{cents}. Now suppose $F = \langle t, u \rangle \cong 2^2$. Then $C_G(t)^0 = A_1 A_5$, $F_4$ or $C_4$ and therefore $C_G(F)^0 = A_1A_3$, $A_1C_3$, $B_4$ or $C_2^2$. The $A_1A_3$ case is ruled out by Lemma \ref{e6order4}. The $B_4$ and $C_2^2$ subgroups are both contained in $D_5$-parabolic subgroups. Therefore $C_G(F)^0 = A_1 C_3$. Finally, if $F$ has a further involution $v$ then $C_G(F)^0 = A_1^2 C_2$, which by Proposition \ref{classirred} is $A_1 A_5$-reducible. \hal

We now let $J$ be the Fitting subgroup of $F$. Since $F$ is a $\{2,3\}$-group, $J$ is non-trivial. 

\begin{lem}
If $J$ is not a $2$-group or a $3$-group, then $F$ is as in Table $\ref{aute6}$. 
\end{lem}

\pf
Under the assumptions of the lemma, $J$ has an element $x$ of order $6$. Then $C_G(x)^0 = A_2A_2$ by Proposition \ref{aute6cent}. If $F \neq \langle x \rangle$ then there exists an element $t \in F \setminus J$ inverting $x$ with $t^2 \in J$. Since $F$ has no element of order $4$ we have $t^2=1$ and $\langle J, t \rangle \cong Dih_{12}$. The element $t$ induces a graph automorphism on both $A_2$ factors and so $C_G(J,t)^0 = A_1 A_1$. Since the two factors are non-conjugate, there is no element swapping them and hence $F = \langle J, t \rangle$. \hal

\begin{lem}
If $J$ is a $2$-group then $F = J$.  
\end{lem}

\pf The possibillities for the 2-group $J$ are in Table \ref{aute6}, from which we see that no element of order $3$ can act as a graph automorphism on $C_G(J)^0$. \hal

\begin{lem}
If $J$ is a $3$-group then $F$ is as in in Table $\ref{aute6}$. 
\end{lem}

\pf 
Suppose $J = \langle x \rangle \cong Z_3$. If $F = J$ then $C_G(F)^0 = A_2^3$ by Proposition \ref{cents}. Otherwise, $F = \langle J, t \rangle \cong Dih_6$ where $t$ is an involution in $N_G(A_2^3) / A_2^3 \cong 2 \times S_3$ by Proposition \ref{normalizers}. This gives two possibilities for $t$. If $t$ is the central involution then $t$ induces a graph automorphism on each factor $A_2$ and $C_G(J,t)^0 = A_1^3$. If $t$ is not central then $C_G(J,t)^0 = A_1 A_1$. 

Now suppose $|J| > 3$ and let $x,y \in J$ with $\la x,y\ra \cong 3^2$. Then $C_G(x) = A_2^3 . 3$ and $y$ cyclically  permutes the three $A_2$ factors, so $C_G(x,y)^0$ is a diagonal subgroup $A_2$. However, the elements in class $3A$ have trace $-3$ on $L(G)$ and so 
\[
\dim C_{L(G)}(S) = \frac{1}{9}(78-8\cdot 3) = 6,
\]
 a contradiction. \hal 

\vspace{4mm}

Finally, we need to prove that all the subgroups $C_G(F)^0$ in Table \ref{aute6} are $G'$-irreducible.

\begin{lem}\label{existe6}
Let $F$ and $C_G(F)^0$ be as in Table $\ref{e6}$. Then $C_G(F)^0$ is $G$-irreducible.
\end{lem}

\pf
This is proved in a similar fashion to Lemma \ref{exist} for most of the subgroups. Specifically, all of the subgroups have no trivial composition factors on $L(G)$ except for $A_1 A_5$, $\bar{A}_1 C_3$ and $\bar{A}_1 A_3$ when $p=3$, all of which have exactly one trivial composition factor. There are no Levi subgroups of $G$ containing a subgroup of type $A_1 A_5$ or $A_1 C_3$. Hence both are $G$-irreducible. 

Now consider $X = A_1 A_3$. Assume $X$ is $G$-reducible and choose a minimal parabolic subgroup $P$ containing $X$. By \cite[Theorem 1]{LSred}, $X$ is contained in a Levi subgroup $L$ of $P$ and by minimality $X$ is $L$-irreducible. Hence $L = D_5 T_1$ or $A_1 A_3 T_2$, where $T_i$ denotes a central torus of rank $i$. The second possibility is ruled out since $X$ has only one trivial composition factor on $L(G)$. So $X$ is an irreducible subgroup of $L' = D_5$. The $A_1$ factor of $X$ is generated by root groups of $D_5$ and so $C_{D_5}(A_1)^0 = A_1 A_3$. Thus $C_{D_5}(X)^0$ contains a subgroup $A_1$, contradicting the $L$-irreducibility of $X$. Hence $X$ is $G$-irreducible, as required.       
\hal

This completes the proof of Theorem \ref{excep} for $G = {\rm Aut}\,E_6$.

\subsection{The case $G=F_4$}

Let $G=F_4$, and let $F$ be a finite subgroup of $G$ such that $C_G(F)^0$ is $G$-irreducible. 

\begin{lem}
If $F$ is an elementary abelian $2$-group, then $F$ and $C_G(F)^0$ are given in Table $\ref{f4}$. 
\end{lem}

\pf
If $F \cong Z_2$ then $C_G(F) = B_4$ or $A_1 C_3$ by Proposition \ref{cents}. Now suppose $F = \langle t, u \rangle \cong 2^2$. Then $u \in C_G(t) = B_4$ or $A_1 C_3$. Therefore $C_G(F)^0 = D_4$ or $A_1^2 B_2$. Now suppose $|F| > 4$. A $2^2$ subgroup of $F$ must contain a $2A$ involution, say $t$, with centralizer $B_4$. Then $B_4 /\langle t \rangle \cong SO_9$ and the image of $F$ in $SO_9$ is elementary abelian by Proposition \ref{spincents}. Since $C_G(F)^0$ is $G$-irreducible, it follows that the image is $\langle u, v \rangle  \cong 2^2$ with $u = (-1^8,1)$ and $v = (-1^4,1^5)$. Therefore $C_{SO_9}(F)^0 = SO_4 SO_4$, and so $C_G(F)^0 = A_1^4$. \hal

\begin{lem}
If $F$ is a $2$-group containing an element $x$ of order $4$, then $F$ and $C_G(F)^0$ are given in Table $\ref{f4}$. 
\end{lem}

\pf
By Proposition \ref{cents}, $C_G(x) = A_1 A_3$. Now suppose $|F| > 4$. Since $A_1 A_3$ is not contained in $A_1 C_3$ it follows that $C_G(x^2) = B_4$. Therefore $F / \langle x^2 \rangle$ is elementary abelian in $SO_9$. Diagonalising as before we may assume $x = (1^3,-1^6)$. Since $C_G(F)^0$ is $G$-irreducible, it must be the case that $|F / \langle x^2 \rangle| = 4$ and a further involution in $F$ is either $u_1 = (-1^6,1^3)$ or $u_2 = (1^5,-1^4)$. In the first case the order of $x u_1$ is $4$ and hence $F \cong Q_8$ with $C_G(F)^0 = B_1^3$. In the second case the order of $x u_2$ is $2$ and hence $F \cong Dih_8$ with $C_G(F)^0 = B_1 B_2$.   
\hal

Now let $J$ be the Fitting subgroup of $F$. Since $F$ has no element of order $6$ it follows that $J$ is either a $2$-group or a $3$-group. 

\begin{lem}
If $J$ is a $2$-group then $F$ and $C_G(F)^0$ are given in Table $\ref{f4}$. 
\end{lem}

\pf 
By the previous two lemmas we may assume that $F$ is not a $2$-group, hence contains an element $x$ of order $3$. The only possibilities for $J$ are $2^2$, $2^3$ or $Q_8$, with $C_G(J)^0 = D_4$,  $A_1^4$ or $B_1^3$, respectively. The last two cases are ruled out since any proper diagonal connected subgroup of $A_1^4$ or $B_1^3$ such that each projection involves no nontrivial field automorphisms is not $B_4$-irreducible by Proposition \ref{classirred}. 

Hence $J \cong 2^2$ and $C_G(J)^0 = D_4$. If $F = \la J,x\ra \cong Alt_4$, then $C_G(F)^0 = A_2$ or $G_2$; and $G_2$ is not possible since it is contained in a Levi subgroup of type $B_3$. And if $F \ne \la J,x\ra$ then $F \cong Sym_4$ and $C_G(F)^0 = A_1$.  \hal

\begin{lem}
If $J$ is a $3$-group then $F$ and $C_G(F)^0$ are given in Table $\ref{f4}$. 
\end{lem}

\pf
Let $J = \langle x \rangle$, so $C_G(x) = A_2 A_2$. Proposition \ref{normalizers} gives   
$N_G(A_2 A_2) / A_2 A_2  = \la t\ra \cong Z_2$, where $t$ acts as a graph automorphism on each factor. Therefore $F = \langle x, t \rangle \cong Dih_6$ with $C_G(F)^0 = A_1 A_1$. 
\hal

As before, the fact that all the subgroups $C_G(F)^0$ in Table \ref{f4} are $G$-irreducible is proved in similar fashion to Lemma \ref{exist}; in particular they all have no trivial composition factors on $L(G)$. 

This completes the proof of Theorem \ref{excep} for $G = F_4$.

\subsection{The case $G = G_2$}

\begin{lem}
Let $F$ be a finite subgroup of $G=G_2$ such that $C_G(F)^0$ is $G$-irreducible. Then $F$ and $C_G(F)^0$ are as in Table $\ref{g2}$. 
\end{lem}

\pf
By Proposition \ref{cents}, non-identity elements of $F$ have order 2 or 3. If $F$ is a $2$-group then $F$ contains an involution $t$ with $C_G(t) = A_1 A_1$; the centralizer of an involution in $A_1 A_1$ is reducible and therefore $F =\la t\ra$. Similarly, if $F$ is a $3$-group then $F = \langle u \rangle \cong 3$ and $C_G(F) = A_2$. The only remaining possibility is $F = \langle t, u \rangle \cong Dih_6$. Since $N_G(A_2) / A_2 \cong 2$, such an example exists and $C_G(F)^0 = A_1$.  
\hal

This completes the proof of Theorem \ref{excep}.

\section{Proof of Proposition \ref{class}} \label{classproof}

In this section we prove the following generalisation of Proposition \ref{class}.

\begin{prop}\label{classgen}
Let $G$ be a classical simple adjoint algebraic group in characteristic $p \ge 0$ with natural module $V$, and let $H = {\rm Aut}\,G$. 
Suppose $F$ is a finite subgroup of $H$ such that $C_G(F)^0$ is $G$-irreducible. Then $F$ is an elementary abelian $2$-group (or a group of order $3$ or $6$ in the case where $G = D_4$), and one of the following holds.
\begin{itemize}
\item[{\rm (i)}] $G = PSL_n$, $F\cap G = 1$, $|F| = 2$ and
\begin{itemize}
\item[ ] if $n$ is even, then $C_G(F)^0 = PSp_n$ or $PSO_n\,(p\ne 2)$;
\item[ ] if $n$ is odd, then $p\ne 2$ and $C_G(F)^0 = PSO_n$.
\end{itemize}
\item[{\rm (ii)}] $G = H = PSp_{2n}$, $p\ne 2$, and  taking preimages in $Sp_{2n}$,  
\[
C_{Sp_{2n}}(F)^0 = \prod_i Sp_{2n_i} = \prod_i Sp(W_i),
\]
where $\sum n_i = n$ and $W_i$ are the distinct weight spaces of $F$ on $V$.
\item[{\rm (iii)}] $G = PSO_n\,(n\ne 8)$, $H = PO_n$, $p\ne 2$, and taking preimages in $O_n$,  
\[
C_{O_n}(F)^0 = \prod_i SO_{n_i} = \prod_i SO(W_i),
\]
where $n_i\ge 3$ for all $i$, $\sum n_i = n$ or $n-1$, and $W_i$ are weight spaces of $F$.
\item[{\rm (iv)}] $G = PSO_{2n}\,(n\ne 4)$, $p=2$, $F\cap G = 1$, $|F| = 2$ and $C_G(F)^0 = SO_{2n-1}$.
\item[{\rm (v)}] $G=D_4=PSO_8$, $H = D_4.Sym_3$, and $F$, $C_G(F)^0$ are as in Table $\ref{d4tab}$.
\end{itemize}
\end{prop}

\pf First suppose $G = PSL_n$. If $F \cap G \ne 1$ then $C_G(F\cap G)^0$ is reducible, so $F\cap G = 1$. Hence $|F|=2$ and now the conclusion in part (i) follows from Proposition \ref{outeraut}. Similarly, if $G = D_n = PSO_{2n}$ with $n\ne 4$ and $p=2$ (so that $H = G.2$), then $F\cap G = 1$, $|F|=2$ and $C_G(F)^0 = B_{n-1}$ by Proposition \ref{outeraut}, as in (iv).

Now suppose $G = PSp_{2n}$. Then $H=G$ and the centralizer in $G$ of any element of order greater than 2 is reducible. Hence $F$ is an elementary abelian 2-group and $p\ne 2$. The preimage $\hat F$ of $F$ in $Sp_n$ must also be elementary abelian, and if we let $W_i\,(1\le i\le k)$ be the weight spaces of $\hat F$ on $V$, then $V = W_1\perp \cdots \perp W_k$ and $C_{Sp_{2n}}(\hat F) = \prod Sp(W_i)$, as in conclusion (ii).

A similar proof applies when $G = PSO_n$ with $n\ne 8$ and $p\ne 2$, giving (iii).

It remains to handle $G = D_4 = PSO_8$. Here $H = G.Sym_3$. If $F \le PO_8 = G.2$ then the above proof shows that $F = 2$ or $2^2$ is as in Table \ref{d4tab}. Now suppose 3 divides $|F|$, so that $F$ contains an element $x$ of order 3 inducing a triality automorphism on $G$. By Proposition \ref{outeraut}, $C_G(x) = G_2$ or $A_2$, with $p\ne 3$ in the latter case. 

If there is an element $y \in C_F(x)\setminus \la x\ra$, then $y \in G_2$ or $A_2$ has irreducible centralizer, which forces $y$ to be an involution in $G_2$. So in this case $F = \la x,y \ra \cong Z_6$ and $C_G(F) = C_{G_2}(y) = \bar A_1A_1$, as in Table \ref{d4tab}. This subgroup has composition factors of dimensions $1$, $3$ and $4$ on $V$, so is $G$-irreducible.

We may now suppose that $C_F(x) = \la x\ra$ and $F \ne \la x\ra$. This implies that $F = \la x,t\ra \cong Dih_6$. If 
$C_G(x) = G_2$ then $t$ must centralize $G_2$, so that $C_G(F)=G_2$. And if $C_G(x) = A_2$ then $t$ induces a graph automorphism on $A_2$ (see Proposition \ref{outeraut}), so $C_G(F) = A_1$ and $p\ne 2$, as in Table \ref{d4tab}. This completes the proof. \hal

\begin{table}
\caption{$G={\rm Aut}\,D_4$:  finite subgroups $F$ with irreducible centralizer} \label{d4tab}
\[
\begin{array}{|c|c|c|}
\hline
F &F\cap G^0 & C_G(F)^0  \\
\hline
2 & 2 &  A_1^4 \,(p\ne 2) \\
  & 1 & B_3 \\
  & 1 & B_1B_2\,(p\ne 2) \\
2^2 & 2 & A_1^2B_1\,(p\ne 2) \\
3 & 1 & G_2 \\
  &  1 & A_2\,(p\ne 3) \\
6 & 2 & \bar A_1A_1\,(p\ne 2) \\
Dih_{6} & 1 & G_2 \\
        && A_1\,(p\ne 2,3) \\
\hline
\end{array}
\]
\end{table}

\section{Tables of results}\label{tabs}

This section consists of the tables referred to in Theorem \ref{excep}.

\begin{table}
\caption{$G=E_8$:  finite subgroups $F$ with irreducible centralizer} \label{e8}
\[
\begin{array}{|c|c|c|c|}
\hline
F & C_G(F)^0 & \hbox{elements of }F & \hbox{ maximal rank} \\
& & & \hbox{overgp. of }C_G(F)^0 \\
\hline
 2 & A_1E_7 & 2A & \\
  & D_8 & 2B &\\
2^2 & D_4^2 & 2B^3 & \\
       & A_1^2D_6 & 2A^2,2B &  \\
2^3 & A_1^4D_4 & 2A^4,2B^3 &  \\
     & A_1^8 & 2B^7 & \\
2^4 & A_1^8 & 2A^8,2B^7 &  \\
4 & A_1A_7 & 2A,4A^2 & \\
 & A_3D_5 & 2B,4B^2 & \\
4\times 2 & A_1^2A_3^2 & 2A^2,2B,4B^4 & \\
Dih_8 & B_1^2B_4 & 2A^4,2B,4B^2 & D_8 \\
         & B_1^2B_2^2 & 2B^5,4B^2 & D_8 \\
          & B_1B_2B_3 & 2A^2,2B^3,4B^2 & D_8 \\
Q_8 & \bar{A}_1D_4 & 2A,4A^6 & A_1A_7 \\
        & B_2^3 & 2B,4B^6 & D_8 \\
      & B_1^3B_3 & 2B,4B^6 & D_8 \\
Dih_8\times 2 & \bar{A}_1^2B_1^2B_2 & 2A^6,2B^5,4B^4 & D_8 \\
4\circ Dih_8 & A_3B_1^3 & 2A^6,2B,4B^8 & D_8 \\
Q_8\times 2 & \bar{A}_1^2B_1^4 & 2A^2,2B,4B^{12} & D_8 \\
2_-^{1+4} & B_1^5 & 2A^{10},2B,4B^{20} & D_8 \\
\hline
3 & A_2E_6 & 3B^2 & \\
   &  A_8 & 3A^2 & \\
3^2 & A_2^4 & 3B^8 & \\
5 & A_4^2 & 5A & \\
\hline
6 & A_1A_2A_5 & 2A, 3B^2, 6A^2 & \\
Dih_6 & A_1F_4 & 2A^3, 3B^2 & A_2E_6 \\
          & A_1C_4 & 2B^3, 3B^2 & A_2E_6 \\
          & B_4 & 2B^3, 3A^2 & A_8 \\
Dih_{12} & \bar{A}_1A_1C_3 & 2A^4, 2B^3, 3B^2, 6A^2 &  A_1A_2A_5 \\
G_{12} & \bar{A}_1A_1A_3 & 2A, 3B^2, 4A^6, 6A^2 & A_1A_2A_5 \\
Alt_4 & A_2G_2 & 2B^3, 3B^8 & D_4^2 \\
    & A_2A_2 & 2B^3, 3A^8 & D_4^2 \\
Sym_4 &  A_1G_2 & 2A^6, 2B^3, 3B^8, 4B^6 &  D_4^2 \\
    &  A_1A_1 & 2B^9, 3A^8, 4B^6 &  D_4^2 \\
Alt_4\times 2 & \bar{A}_1A_1A_2 & 2A^4, 2B^3, 3B^8, 6A^8 & A_1^4D_4 \\
Sym_4\times 2 & \bar{A}_1A_1A_1 & 2A^{12}, 2B^7, 3B^8, 4B^{12}, 6A^8 & A_1^4D_4 \\
SL_2(3) & \bar{A}_1A_2 & 2A, 3B^8, 4A^6, 6A^8 & A_1A_7 \\
3^2.2 & A_1^4 & 2B^9, 3B^8 & A_2^4  \\
     & A_1A_2A_2 & 2A^3, 3B^8, 6A^6 & A_2^4 \\
3^2.4 & A_1^2 & 2B^9, 3B^8, 4B^{18} & A_2^4 \\
3^2.2^2 & A_1^2A_1 & 2A^6, 2B^9, 3B^8, 6A^{12} & A_2^4 \\
3^2.Dih_8 & A_1^2 & 2A^{12}, 2B^9, 3B^8, 4B^{18}, 6A^{24} & A_2^4 \\
\hline
Dih_{10} & B_2^2 & 2B^5, 5A^4 & A_4^2 \\
Frob_{20} & B_2 & 2B^5, 4B^{10}, 5A^4 & A_4^2 \\
\hline
Alt_5 & A_1 & 2B^{15}, 3B^{20}, 5A^{24} & A_4^2 \\
Sym_5 & A_1 & 2A^{10}, 2B^{15}, 3B^{20}, 4B^{30}, 5A^{24}, 6A^{20} & A_4^2 \\
\hline
\end{array}
\]
\end{table}

\begin{table}
\caption{$G=E_7$ (adjoint):  finite subgroups $F$ with irreducible centralizer} \label{e7}
\[
\begin{array}{|c|c|c|c|}
\hline
F & C_G(F)^0 & \hbox{elements of }F & \hbox{ maximal rank} \\
& & & \hbox{overgp. of }C_G(F)^0 \\
\hline
2 & A_1D_6 & 2A & \\
 & A_7 & 2B & \\
2^2 & A_1^3D_4 & 2A^3 & \\
     & D_4 & 2B^3 & A_7 \\
2^3 & A_1^7 & 2A^7 &  \\
4 & A_1A_3^2 & 2A,4A^2 & \\
Dih_8 & \bar{A}_1B_1^2B_2 & 2A^5,4A^2 & A_1D_6 \\
Q_8 & \bar{A}_1B_1^4 & 2A,4A^6 & A_1D_6 \\
\hline
3 & A_2A_5 & 3A & \\
Dih_6 & A_1C_3 & 2A^3, 3A^2 & A_2A_5 \\
   & A_1A_3 & 2B^3, 3A^2 & A_2A_5 \\
Alt_4 & A_1A_2 & 2A^3,3A^8 & A_1^3D_4 \\
   & A_2 & 2B^3, 3A^8 & A_7 \\
Sym_4 & A_1A_1 & 2A^9, 3A^8, 4A^6 & A_1^3D_4 \\
\hline
\end{array}
\]
\end{table}

\begin{table}
\caption{$G=E_6$:  finite subgroups $F$ with irreducible centralizer} \label{e6}
\[
\begin{array}{|c|c|c|c|}
\hline
F & C_G(F)^0 & \hbox{elements of }F & \hbox{ maximal rank} \\
& & &  \hbox{overgp. of }C_G(F)^0 \\
\hline
2 &  A_1A_5 & 2A & \\
3 &  A_2^3 & 3A^2 & \\
Dih_6 &  A_1A_1 & 2A^3, 3A^2 & A_2^3 \\
\hline
\end{array}
\]
\end{table}

\begin{table}
\caption{$G={\rm Aut}\,E_6$:  finite subgroups $F$ with irreducible centralizer} \label{aute6}
\[
\begin{array}{|c|c|c|c|c|}
\hline
F &F\cap G^0 & C_G(F)^0 & \hbox{elements of }F & \hbox{ maximal rank} \\
& & & &  \hbox{overgp. of }C_G(F)^0 \\
\hline
2 & 2 & A_1A_5 & 2A & \\
  & 1 & F_4 & 2B & \\
  & 1 & C_4 & 2C & \\
4 & 2 & \bar{A}_1A_3 & 2A, 4A^2 & A_1A_5 \\
2^2 & 2 & \bar{A}_1C_3 & 2A, 2B, 2C & A_1A_5 \\
3 & 3 & A_2^3 & 3A^2 & \\
Dih_6 & Dih_6 & A_1A_1 & 2A^3, 3A^2 & A_2^3 \\
Dih_{12} & Dih_6 & A_1A_1 & 2A^3, 2B, 2C^3, 3A^2, 6A^2 & A_2^3 \\
\hline
\end{array}
\]
\end{table}

\begin{table}
\caption{$G=F_4$:  finite subgroups $F$ with irreducible centralizer} \label{f4}
\[
\begin{array}{|c|c|c|c|}
\hline
F & C_G(F)^0 & \hbox{elements of }F & \hbox{ maximal rank} \\
& & &  \hbox{overgp. of }C_G(F)^0 \\
\hline
2 & B_4 & 2A & \\
   & A_1C_3 & 2B & \\
2^2 & A_1^2C_2 & 2A, 2B^2 & \\
      & D_4 & 2A^3 & \\
2^3 & A_1^4 & 2A^3, 2B^4 & \\
4 &  A_1A_3 & 2A, 4A^2 & \\
Dih_8   & B_1B_2 & 2A^3, 2B^2, 4A^2 & D_4 \\
Q_8 & B_1^3 & 2A, 4A^6 & B_4 \\
3 & A_2A_2 &  3A^2 & \\
Dih_6 & A_1A_1 & 2B^3, 3A^2 & A_2A_2 \\
Alt_4 & A_2 & 2A^3, 3A^8 & D_4 \\
Sym_4 & A_1 & 2A^3, 2B^6, 3A^8, 4A^6 & D_4 \\
\hline
\end{array}
\]
\end{table}

\begin{table}
\caption{$G=G_2$:  finite subgroups $F$ with irreducible centralizer} \label{g2}
\[
\begin{array}{|c|c|c|c|}
\hline
F & C_G(F)^0 & \hbox{elements of }F & \hbox{ maximal rank} \\
& & &  \hbox{overgp. of }C_G(F)^0 \\
\hline
2 & A_1A_1 & 2A & \\
3 & A_2 & 3A^2 & \\
Dih_6 & A_1 & 2A^3, 3A^2 & A_2 \\
\hline
\end{array}
\]
\end{table}

\vspace*{0.5cm}

\noindent
M.W. Liebeck, Imperial College, London SW7 2AZ, UK, 

\no m.liebeck@imperial.ac.uk

\vspace{4mm}
\no A.R. Thomas, Heilbronn Institute for Mathematical Research, University of Bristol, Bristol, UK

\no adam.thomas@bristol.ac.uk


\begin{thebibliography}{99}
\bibitem{AS} M.\ Aschbacher and G.M.\ Seitz, Involutions in Chevalley groups 
over fields of even order, 
{\it Nagoya Math. J.}  {\bf 63} (1976), 1--91. 

\bibitem{BT} A. Borel and J. Tits, \'El\'ements 
unipotents et
sousgroupes paraboliques de groupes r\'eductifs, {\it Invent. Math.}
{\bf 12} (1971), 95--104.


\bibitem{bourbaki}
N.~Bourbaki, \emph{Groupes et {A}lg\`ebres de {L}ie \textup{(Chapters 4,5,6)}}, Hermann, Paris, 1968.

\bibitem{Car} R.W. Carter, Conjugacy classes in the Weyl group, {\it Compositio Math.} {\bf 25} (1972), 1--59. 

\bibitem{CG} A.M. Cohen and R.L. Griess, On finite simple subgroups of the complex
Lie group of type $E_8$, {\it Proc. Symp. Pure Math.} {\bf 47} (1987), 367--405.

\bibitem{CW}
A.M. Cohen and D.B. Wales, Finite subgroups of $E_6(\C)$ and $F_4(\C)$,
{\it Proc. London Math. Soc.} {\bf 74} (1997), 105--150.

\bibitem{CLSS} A.M. Cohen, M.W. Liebeck, J. Saxl and G.M. Seitz, The local maximal subgroups
of exceptional groups of Lie type, finite and algebraic, {\it Proc. London
Math. Soc.} {\bf 64} (1992), 21--48.

\bibitem{GLS}
D.\ Gorenstein, R.\ Lyons and R.\ Solomon, \emph{The classification of 
the finite simple groups. Number 3. Part I. Chapter A. 
Almost simple $K$-groups.} Mathematical Surveys and Monographs, 
40.3. American Mathematical Society, Providence, RI, 1998.

\bibitem{kleidman} P.B. Kleidman, {The maximal subgroups of the finite $8$-dimensional orthogonal groups $P\Omega^+_8(q)$ and of their automorphism groups}, {\it J. Algbera} {\bf 110} (1987), 173--242.

\bibitem{geomded} M.W. Liebeck and G.M. Seitz, Maximal subgroups of
exceptional groups of Lie type, finite and algebraic, {\it Geom. Ded.}
{\bf 36} (1990), 353--387.

\bibitem{LSred} M.W. Liebeck and G.M. Seitz, Reductive subgroups of exceptional algebraic groups, Mem. Amer. Math. Soc. {\bf 121} (1996), no. 580.


\bibitem{LSbk} M.W. Liebeck and G.M. Seitz,
{\em Unipotent and nilpotent classes in simple algebraic groups and Lie algebras\/},
Mathematical Surveys and Monographs,  Vol.180, American Math. Soc., Providence, RI, 2012.

\bibitem{LT} M.W. Liebeck and D.M. Testerman, Irreducible subgroups of algebraic groups, 
{\it Quarterly J. Math} {\bf 55} (2004), 47--55.

\bibitem{MP} R.V. Moody and J. Patera, Characters of elements of finite order in Lie groups, {\it SIAM J. Algebra Discrete Math.} {\bf 5} (1984), 359--383. 

\bibitem{ser04}
J.-P. Serre, \emph{Compl\`{e}te r\'{e}ductibilit\'{e}}, Ast\'{e}risque \textbf{299} (2005), Exp. No. 932, S\'{e}minaire Bourbaki. Vol. 2003-2004.

\bibitem{Tho} A.R. Thomas, The irreducible subgroups of exceptional algebraic groups, preprint. 


\end{thebibliography}
\end{document}